\pdfoutput=1
\RequirePackage{ifpdf}
\ifpdf 
\documentclass[pdftex]{sigma}
\else
\documentclass{sigma}
\fi

\usepackage{tikz}
\usetikzlibrary{decorations.pathreplacing,arrows,positioning,matrix}
\tikzset{>=stealth}

\numberwithin{equation}{section}
\newtheorem{Theorem}{Theorem}[section]
\newtheorem{Corollary}[Theorem]{Corollary}
\newtheorem{Lemma}[Theorem]{Lemma}
{\theoremstyle{definition}
\newtheorem{Definition}[Theorem]{Definition} }

\newcommand{\pair}[1]{\ensuremath{\left\langle #1 \right\rangle}}

\DeclareMathOperator{\Au}{\mathbf{A}}
\DeclareMathOperator{\Wu}{\mathbf{W}}
\DeclareMathOperator{\End}{End}

\DeclareMathOperator{\Hom}{Hom}

\definecolor{rd}{rgb}{0.0, 0.0, 0.9}
\definecolor{io}{rgb}{1.0, 0.65, 0.0}
\definecolor{lb}{rgb}{0.48, 0.45, 0.8}

\usepackage{ifthen}
\newcommand{\Z}[3]{\ifthenelse{#1<#2}{\ensuremath{Z^{#3}_{#1,#2}}}{\ensuremath{\underline{Z^{#3}_{#1,#2}}}}}

\begin{document}

\newcommand{\arXivNumber}{1410.7593}

\allowdisplaybreaks

\renewcommand{\PaperNumber}{053}

\FirstPageHeading

\ShortArticleName{Constructing Involutive Tableaux with Guillemin Normal Form}
\ArticleName{Constructing Involutive Tableaux\\ with Guillemin Normal Form}

\Author{Abraham D.~SMITH}
\AuthorNameForHeading{A.D.~Smith}

\Address{Department of Mathematics, Statistics and Computer Science, University of Wisconsin-Stout,\\ Menomonie, WI  54751-2506, USA}
\Email{\href{mailto:adsmith@member.ams.org}{adsmith@member.ams.org}}
\URLaddress{\url{http://www.curieux.us/abe/}}

\ArticleDates{Received December 15, 2014, in f\/inal form July 01, 2015; Published online July 09, 2015}

\Abstract{Involutivity is the algebraic property that guarantees solutions to an
analytic and torsion-free exterior dif\/ferential system or partial dif\/ferential
equation via the Cartan--K\"ahler theorem.  Guillemin normal form establishes
that the prolonged symbol of an involutive system admits a commutativity
property on certain subspaces of the prolonged tableau.  This article examines
Guillemin normal form in detail, aiming at a more systematic approach to
classifying involutive systems.
The main result is an explicit quadratic condition for involutivity of the type
suggested but not completed in Chapter~IV, \S~5 of the book \emph{Exterior
Differential Systems} by Bryant, Chern, Gardner, Goldschmidt, and Grif\/f\/iths.
This condition enhances Guillemin normal form and characterizes involutive
tableaux.}

\Classification{58A15; 58H10}

\section{Background}
Fix real or complex vector spaces (or bundles, etc.) $W$ and $V$ of dimension $r$ and $n$,
respectively. A
\textit{tableau}\footnote{The literature contains many
distinct-yet-related meanings of \emph{tableau} and \emph{symbol}.
Usage here is consistent with~\cite{BCGGG}.  In~\cite{Guillemin1968}, Guillemin
uses the phrase ``a subspace of~$\Hom(V,W)$'' where we say \emph{tableau}.}
is a vector space~$A$ with an exact sequence
\begin{gather*}
0 \longrightarrow A \longrightarrow W \otimes V^* \overset{\sigma}{\longrightarrow} H^{1}(A) \longrightarrow 0.
\label{spencer1}
\end{gather*}
A homomorphism $\sigma$ with kernel $A$ is called a \textit{symbol}, taking
values in the cokernel $H^{1}(A)=A^\perp$.

Given bases of $(w_a)$ of $W$ and $(u^i)$ of $V^*$, elements of $A$ are expressed as $r \times n$
matrices $\pi = \pi^a_i (w_a \otimes u^i)$.
In $(\pi^a_i)$, the independent generators of $A$ appear in some $s_1$ entries of the f\/irst
column, some $s_2$ entries of the second column, and so on to some $s_n$
entries in the last column.
These numbers are constant over a Zariski-open subset of the bases (called the
\textit{generic} bases) of $V^*$, and those constant values satisfy $s_1
\geq s_2 \geq s_3 \geq \cdots \geq s_n \geq 0$.  That is, for a generic basis
of $V^*$, the independent generators of $A$ are packed to the left in $(\pi^a_i)$.
Similarly, there is a Zariski-open subset of generic bases $(w_a)$ of $W$
wherein the independent
generators of $A$ are packed to the top.  If both $(u^i)$ and $(w_a)$ are
generic, then the generators of $A$ appear in the f\/irst $s_1$ entries of the
f\/irst column, the f\/irst $s_2$ entries of the second column, and so on.
Hence, $\dim A = s_1 + s_2 + \cdots + s_n$.

Given a tableau $A$, consider the tableau $A^{(1)} \to  A\otimes V^*$ given as the kernel of
the map $\delta_\sigma\colon A \otimes V^* \to W \otimes \wedge^2 (V^*)$ by
$\delta_\sigma = (1_W \otimes \delta)\circ(\sigma \otimes 1_{V^*})$ where $\delta$ is
the skewing map $\delta\colon V^* \otimes V^* \to V^* \wedge V^*$.
This $A^{(1)}$ is called the (f\/irst) \textit{prolonged tableau}, and the map
$\delta_\sigma$ is the (f\/irst) \textit{prolonged symbol}.  The cokernel is
written $H^{2}(A)$,
\begin{gather*}
0 \longrightarrow A^{(1)} \longrightarrow A \otimes V^* \overset{\delta_\sigma}{\longrightarrow} W \otimes
\wedge^2 V^* \longrightarrow H^2(A) \longrightarrow 0.
\label{spencer2}
\end{gather*}
It is a standard exercise to show
$\dim A^{(1)} \leq s_1 + 2s_2 + \cdots + n s_n$.
The most interesting property to study for a tableau is
\begin{Definition}[involutivity]
A tableau is called \emph{involutive} if and only if equality holds in the relation
$\dim A^{(1)} \leq s_1 + 2s_2 + \cdots + ns_n$.
\end{Definition}

Verif\/ication of this equality is called \textit{Cartan's test}. Its most
frequent application is to check whether an overdetermined system of partial
dif\/ferential equations (or an exterior dif\/ferential system) admits solutions to
the analytic Cauchy initial-value problem. Involutivity is a coor\-di\-na\-te-invariant property of a tableau,
and techniques from homological algebra comprise most of the
recent literature, for example \cite{CGG2009, Kruglikov2007, Malgrange2003}.

In \cite{Guillemin1968}, Guillemin showed involutive tableaux admit
a partially commuting ``normal form'' for the symbol map, using Quillen's results on the
exactness of Spencer cohomology from \cite{Quillen1964}.  Guillemin normal form
was reconsidered for exterior dif\/ferential systems in \cite{Yang1987} and
\cite{BCGGG}.  This
article\footnote{This article began as an appendix to \cite{Smith2014b}, so there is some overlap in the presentation.
However, be aware that some indices~-- such as~$i$,~$j$~-- and some notations~-- such
as~$Y^*$~-- dif\/fer between the two articles.}
is intended to extend and clarify Guillemin normal form
to aide future computational and theoretical applications regarding the
geometry of partial
dif\/ferential equations.  The main result is Theorem~\ref{thm:gnf+}.

\section{Endovolutivity}
\label{sec:notation}
Let $\ell$ denote the index of the last non-zero Cartan character, $s_\ell$.
Permanently reserve the index ranges
$i,j,k,l \in \{ 1, \ldots, n\}$ and
$\lambda,\mu \in \{ 1, \ldots, \ell\}$ and
$\varrho, \varsigma \in \{ \ell+1, \ldots, n\}$ and
$a,b,c,d \in \{1, \ldots, r\}$.

Let $(u^i)$ and $(w_a)$ denote generic bases of $V^*$ and $W$, respectively, so that an
element $\pi \in W \otimes V^*$ is written as a matrix $\pi = \pi^a_i
(w_a \otimes u^i)$.
Let $U^*$ denote the $\ell$-dimensional subspace spanned by $u^1, \ldots,
u^\ell$. Let $Y^*$ denote the complementary subspace spanned by $u^{\ell+1},
\ldots, u^n$.
We also denote the basis-dual spaces $Y=(U^*)^\perp=\pair{u_{\ell+1},\ldots,u_n}$ and $U = (Y^*)^\perp = \pair{u_1,
\ldots, u_\ell}$.

\begin{figure}
\centering
\begin{tikzpicture} [scale=0.5]
\fill[io] (0,0) -- (10, 0) -- (10,-2) -- ( 9,-2) -- ( 9,-2) -- ( 8,-2) -- ( 8,-3) -- ( 7,-3) -- ( 7,-4) -- ( 6,-4) -- ( 6,-4) -- ( 5,-4) -- ( 5,-6) -- ( 4,-6) -- ( 4,-6) -- ( 3,-6) -- ( 3,-8) -- ( 2,-8) -- ( 2,-10) -- ( 1,-10) -- ( 1,-10) -- ( 0,-10) -- cycle;
\draw[very thick]     (0,0) -- (0,-12) -- (14,-12) -- (14,0) -- cycle;
\draw[dotted]     (10,0) -- (10,-12);
\draw (0,-2) node [left=1pt,black] {$s_\ell$};
\draw (0,-10) node [left=1pt,black] {$s_1$};
\draw (0,-6) node [left=1pt,black] {$s_\lambda$};
\draw (0,-3) node [left=1pt,black] {$s_i$};
\draw (0,0) node [above=1pt,black] {$1$};
\draw (3.5,0) node [above=1pt,black] {$\lambda$};
\draw (7.5,0) node [above=1pt,black] {$i$};
\draw (10,0) node [above=1pt,black] {$\ell$};
\draw (14,0) node [above=1pt,black] {$n$};
\draw  (3.5,-4) node {\Large $\pi^b_\lambda$};
\draw  (7.5,-7) node {\Large $\pi^a_i$};
\draw[very thick, ->,left] (7.0,-7) to node {\small $B^{a,\lambda}_{i,b}~$} (4,-4.3);
\draw[very thick, ->] (7.0,-7) -- (5,-1.5);
\draw[very thick, ->] (7.0,-7) -- (0.5,-7.5);
\draw[very thick, ->] (7.0,-7) -- (1.5,-8.5);
\draw[very thick, ->] (7.0,-7) -- (0.5,-9.5);
\draw[very thick, ->] (7.0,-7) -- (6.5,-2.5);
\end{tikzpicture}
\caption{A tableau in coordinates, with Cartan characters $s_1 \geq s_2 \geq
\cdots \geq s_\ell$.  The upper-left shaded entries are
independent generators.  The lower-right entries depend on them via $\pi^a_i =
B^{a,\lambda}_{i,b}\pi^b_\lambda$, summed as in~\eqref{eq:symrels}.}\label{figtab}
\end{figure}
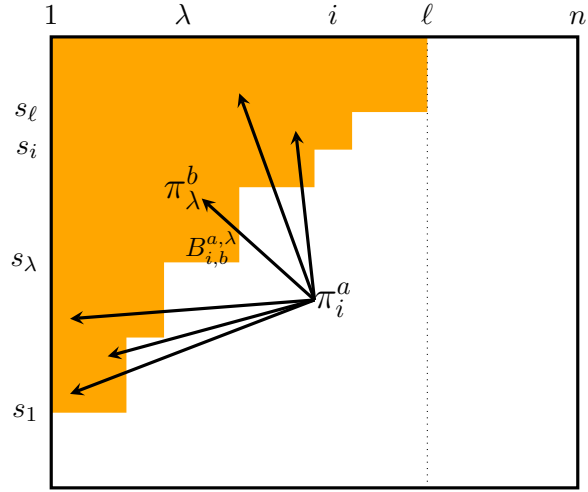

Using generic bases, the symbol $\sigma$ can be expressed as a minimal system of equations of the form
\begin{gather}
\big\{ 0 = \pi^a_i - B^{a,\lambda}_{i,b} \pi^b_\lambda \big\}_{s_i < a},
\label{eq:symrels}
\end{gather}
where $B^{a,\lambda}_{i,b}=0$ unless
$\lambda \leq i$ and  $b \leq s_\lambda$ and $s_i < a $.  See
Fig.~\ref{figtab}.
Using the coef\/f\/icients $B^{a,\lambda}_{i,b}$, we def\/ine an element of
\[
V^* \otimes V \otimes W \otimes W^* \cong \End(V^*) \otimes \End(W)
\] by
the tensorial expression (summed over all $\lambda$, $i$, $b$, as well as~$a$ as
shown)
\begin{gather}
\sum_{a \leq s_i} \delta^\lambda_i \delta^a_b \big(w_a \otimes w^b\big) \otimes \big(u^i
\otimes u_\lambda\big) +
\sum_{a>s_i} B^{a,\lambda}_{i,b} (w_a \otimes w^b) \otimes \big(u^i \otimes
u_\lambda\big).
\label{IB}
\end{gather}
Then, for each $\varphi=\varphi_i u^i \in V^*$, there is a homomorphism $B(\varphi)\colon V \to
\End(W)$ def\/ined by~\eqref{IB} as
\begin{gather*}
B(\varphi)(v) =
\sum_{a \leq s_\lambda}
\varphi_\lambda v^\lambda \delta^a_b \big(w_a \otimes w^b\big) +
\sum_{a>s_i} \varphi_\lambda B^{a,\lambda}_{i,b} v^i \big(w_a \otimes w^b\big),
\label{IBphi}
\end{gather*}
where $v^i = u^i(v)$.  Note that $B(\varphi) = B(\xi)$ if $\varphi_\lambda = \xi_\lambda$ for all
$\lambda$, so~\eqref{IB} is really an element of $V^* \otimes U \otimes
\End(W)$.
We write $B^\lambda_i$ for $B(u^\lambda)(u_i)$, but note that
$(B_\lambda^i)^a_b$ is not
quite the same as~$B^{a,\lambda}_{i,b}$  due to the identity term in~\eqref{IB}; in particular,
\begin{gather*}
\big(B^\lambda_i\big)^a_b =
\begin{cases}
\delta^a_b, & \text{if $\lambda=i$ and $a \leq s_\lambda$},\\
B^{a,\lambda}_{i,b}, & \text{if $a>s_i$}.
\end{cases}
\label{IBcomp}
\end{gather*}

For each $i$, use the bases $(u^i)$ and $(w_a)$ to def\/ine the subspaces
\begin{gather*}
\Wu^-_i = \big\{ z =w_a z^a\colon \; z^a = 0\; \forall\, a > s_i \big\}, \qquad
\Wu^+_i = \big\{ z =w_a z^a\colon \;  z^a = 0\; \forall\, a \leq s_i \big\}.
\end{gather*}
So that $W =  \Wu^-_i \oplus \Wu^+_i$ and $\Wu^-_1 \supset
\Wu^-_2 \supset \cdots \supset \Wu^-_n$ according to $s_1 \geq s_2
\cdots \geq s_n$.  Of course, for $\varrho > \ell$, we have
$\Wu^-_\varrho = 0$.  For each $\lambda$, consider also the subspace
\begin{gather*}
\Au^-_\lambda = \big\{ \pi = B\big(u^\lambda\big)(\cdot)z,\; z \in \Wu^-_\lambda\big\} \subset
A.
\end{gather*}
The symbol relations \eqref{eq:symrels} say that the coef\/f\/icients~$\pi^a_i$ of
$\pi \in \Au^-_\lambda$ are
determined by the choice of $z \in \Wu^-_\lambda$, so
$\Au^-_\lambda$ and $\Wu^-_\lambda$ are isomorphic
via the projection onto f\/irst~$s_\lambda$ entries in column~$\lambda$.
Using this basis and isomorphism, there is a decomposition
\begin{gather}
A = \bigoplus_{\lambda=1}^\ell \Au^-_\lambda \cong \bigoplus_{\lambda=1}^\ell
\Wu^-_\lambda.
\label{Adecomp}
\end{gather}
Specif\/ically, if $\pi = \pi^a_i (w_a \otimes u^i) \in A$, then let
\begin{gather*}
z_\lambda = \sum_a z^a_\lambda w_a \in W \qquad \text{for} \quad z^a_\lambda =
\begin{cases}
\pi^a_\lambda,& a \leq s_\lambda,\\
0, & \text{otherwise}.
\end{cases}
\label{zdecomp}
\end{gather*}
So, the decomposition~\eqref{Adecomp} yields
\begin{gather*}
\pi = \sum_\lambda \pi_\lambda = \sum_\lambda B\big(u^\lambda\big)(\cdot)z_\lambda \in
W \otimes V^*.
\end{gather*}
Since $\dim \Wu^-_\lambda = s_\lambda$, this is a more precise version of the statement
that, for a generic f\/lag, the tableau matrix has~$s_1$ generators in the f\/irst
column, $s_2$ in the second column, and so on until the f\/inal~$s_\ell$
generators in the~$\ell$ column.

For any $\varphi = \varphi_i u^i \in V^*$, use the bases $(u^i)$ and $(w_a)$ to def\/ine the subspaces
\begin{gather*}
\Wu^-(\varphi)  = \Wu^-_{\min\{i\colon  \varphi_i \neq 0\}},\qquad
\Au^-(\varphi)  = \big\{ \pi = B(\varphi)(\cdot)z,\; z \in \Wu^-(\varphi)\big\}.
\end{gather*}

\begin{figure}
\centering
\begin{tikzpicture} [scale=0.5]
\fill[io] (0,0) -- (10, 0) -- (10,-2) -- ( 9,-2) -- ( 9,-2) -- ( 8,-2) -- ( 8,-3) -- ( 7,-3) -- ( 7,-4) -- ( 6,-4) -- ( 6,-4) -- ( 5,-4) -- ( 5,-6) -- ( 4,-6) -- ( 4,-6) -- ( 3,-6) -- ( 3,-8) -- ( 2,-8) -- ( 2,-10) -- ( 1,-10) -- ( 1,-10) -- ( 0,-10) -- cycle;
\draw[very thick]     (0,0) -- (0,-12) -- (14,-12) -- (14,0) -- cycle;
\draw[dotted]     (10,0) -- (10,-12);
\draw (0,-2) node [left=1pt,black] {$s_\ell$};
\draw (0,-10) node [left=1pt,black] {$s_1$};
\draw (0,-6) node [left=1pt,black] {$s_\lambda$};
\draw (0,-3) node [left=1pt,black] {$s_i$};
\draw (0,0) node [above=1pt,black] {$1$};
\draw (3.5,0) node [above=1pt,black] {$\lambda$};
\draw (7.5,0) node [above=1pt,black] {$i$};
\draw (10,0) node [above=1pt,black] {$\ell$};
\draw (14,0) node [above=1pt,black] {$n$};
\draw  (3.5,-3.5) node {\small $\Wu^-_\lambda$};
\draw  (7.5,-2) node {\small $\Wu^-_i$};
\draw  (8.5,-4.5) node[rotate=90] {\small $\Wu^+_i \cap \Wu^-_\lambda$};
\draw  (8.5,-9.5) node[rotate=90] {\small $\Wu^+_i \cap \Wu^+_\lambda$};
\draw[very thick] (7,-6) -- (7,-12) -- (8,-12) -- (8,-6) -- cycle;
\draw[very thick] (3,0) -- (3,-6) -- (4,-6) -- (4,0) -- cycle;
\draw[very thick] (7,0) -- (7,-6) -- (8,-6) -- (8,0) -- cycle;
\draw[dotted] (7,-3) -- (8,-3);
\draw[very thick, ->] (4.0,-3.5) to node[above] {$0$} (6.9,-2);
\draw[very thick, ->] (4.0,-3.5) to node[below] {$B^\lambda_i$ } (6.9,-4.5);
\draw[very thick, ->] (4.0,-3.5) to node[below] {$0$} (6.9,-9);
\end{tikzpicture}
\caption{The map $B^\lambda_i$ for an endovolutive tableau.}
\label{figB}
\end{figure}
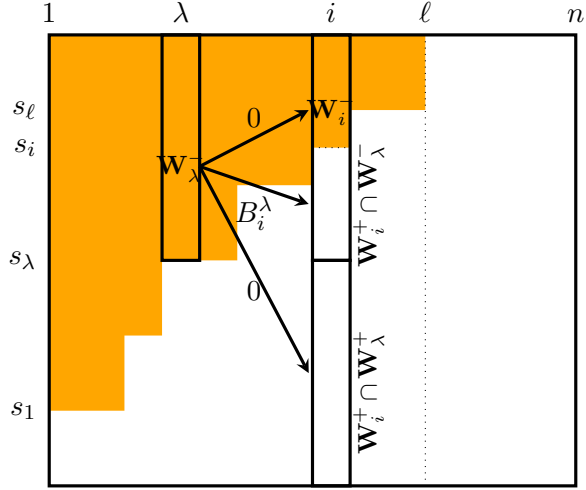

\begin{Definition}[endovolutivity]\label{endodef}
A tableau expressed in generic bases of $V^*$ and $W$ as \eqref{eq:symrels}
is called \emph{endovolutive} if and only if
$B^{a,\lambda}_{i,b}=0$ for all $a> s_\lambda$ using those bases.
See Fig.~\ref{figB}.
\end{Definition}
Endovolutivity is not an invariant property; a tableau could be
endovolutive in generic ba\-ses~$(u^i)$ and~$(w_a)$ but not endovolutive in
other generic bases~$(\tilde{u}^i)$ and~$(\tilde{w}_a)$.  However,
Lemma~\ref{lem:replacement} provides that endovolutivity is invariant under the
Borel subgroup of ${\rm GL}(V^*)$ that preserves the basis f\/lag
\[0 \subset \pair{u^1} \subset
\pair{u^1,u^2} \subset \cdots \subset \pair{u^1, \ldots, u^n} = V^*.\]
\begin{Lemma}\label{lem:replacement}
Suppose a tableau $A$ is endovolutive in generic bases $(u^i)$ and $(w_a)$ for
$V^*$ and $W$.  If $(\tilde{u}^i)$ is another basis of $V^*$ related to $(u^i)$ through a
transformation $\tilde{u}^i = g^i_j u^j$ with $g^i_j=0$ if $i>j$
$($upper-triangular$)$, then the basis $(\tilde{u}^i)$ is also generic for $A$, and
$A$ is endovolutive in $(\tilde{u}^i)$ and $(w_a)$.
\end{Lemma}

The proof is immediate by comparing Figs.~\ref{figtab} and~\ref{figB}, as
the columns of $(\pi^a_i)$ are replaced with linear combinations of columns to
their \emph{right}.

The property of endovolutivity is discussed but not named on page~147 (page~127 in
the online version) of \cite{BCGGG} and in Section~1.2 of~\cite{Yang1987}.
Endovolutivity is so-named here because of Lemma~\ref{lem:endolemma}.

\begin{Lemma}\label{lem:endolemma}
A tableau expressed in generic bases is endovolutive if and only if
$B(\varphi)(v)$ is an endomorphism of $\Wu^-(\varphi)$ for all
$\varphi \in V^*$ and $v \in V$.
\end{Lemma}

\begin{proof}
Suppose $A$ is endovolutive in the given generic bases $(u^i)$ and $(w_a)$.
Fix $\varphi \in V^*$, and let $\kappa = \min
\{ i\colon \varphi_i \neq 0\}$ so that $\Wu^-(\varphi) = \Wu^-_\kappa$.
Replacing $u^\kappa$ with $\varphi$ is an upper-triangular
change-of-basis, so $B(\varphi)(u_j) = \sum\limits_{\lambda \geq \kappa} \varphi_\lambda B^\lambda_j$
preserves $\Wu^-(\varphi)=\Wu^-_{\kappa}$ for all $j$.
Conversely, each~$u^\lambda$ is a particular
choice of $\varphi$, and $\Wu^-(u^\lambda) = \Wu^-_\lambda$, so
$B(u^\lambda)(v)$ is an endomorphism of $\Wu^-(u^\lambda)$ if and only if
$B^{a,\lambda}_{i,b}=0$ for all $a > s_\lambda$.
\end{proof}

When considering endovolutive tableaux, it useful to arrange the symbol
endomorphisms  as an $\ell \times n$ array of $r \times r$ matrices:
\begin{gather}
\vec{B}=
\begin{bmatrix}
I_{s_1}& B^1_2   & B^1_3    & \cdots & B^1_\ell & \cdots & B^1_n \\
0      & I_{s_2} & B^2_3    & \cdots & B^2_\ell & \cdots & B^2_n \\
0      & 0       & I_{s_3}  & \cdots & B^3_\ell & \cdots & B^3_n \\
       &         &           & \ddots &         & B^\lambda_i & \vdots   \\
0      & 0       & 0         & 0      & I_{s_\ell} & \cdots & B^\ell_n
\end{bmatrix}.
\label{bigB}
\end{gather}
In \eqref{bigB}, each $r \times r$ matrix in row $\lambda$ is 0 outside the
upper-left $s_\lambda \times s_\lambda$ part.
For example, the symbol of the endovolutive tableau
\begin{gather*}
\pi=\begin{bmatrix}
\pi^1_1 & \pi^1_2 & \pi^1_3\\
\pi^2_1 & \pi^2_2 & Q_4 \pi^1_2 + Q_5 \pi^2_2 + T_1 \pi^1_1 + T_2 \pi^2_1 + T_3
\pi^3_1\\
\pi^3_1 & P_1\pi^1_1 + P_2 \pi^2_1 + P_3 \pi^3_1 & R_1 \pi^1_1 + R_2 \pi^2_1 + R_3
\pi^3_1
\end{bmatrix},
\end{gather*}
with $(s_1,s_2,s_3) = (3,2,1)$ will be arranged as
\begin{alignat*}{4}
&\begin{pmatrix}
1 & 0 & 0 \\
0 & 1 & 0 \\
0 & 0 & 1
\end{pmatrix},
\qquad &&
\begin{pmatrix}
0 & 0 & 0 \\
0 & 0 & 0 \\
P_1 & P_2 & P_3
\end{pmatrix},
\qquad & &
\begin{pmatrix}
0 & 0 & 0 \\
T_1  & T_2 & T_3 \\
R_1 & R_2 & R_3
\end{pmatrix}, &\\
& &&
\begin{pmatrix}
1 & 0 & ~  \\
0 & 1 & ~  \\
\phantom{T_1}  & \phantom{T_2}  & \phantom{T_3}
\end{pmatrix},
\qquad &&
\begin{pmatrix}
0 & 0 & ~  \\
Q_4 & Q_5 & \phantom{Q_6}  \\
~  & ~  & ~
\end{pmatrix},& \\
&  &&   &&
\begin{pmatrix}
1 & ~  & ~  \\
~  & ~  & ~  \\
\phantom{T_2}  & \phantom{T_2}  & \phantom{T_3}
\end{pmatrix}.&
\end{alignat*}

Note that a change of basis in $V^*$ causes \eqref{bigB} to change by a block-wise
conjugation. The result is not guaranteed to be endovolutive unless the change of basis
is upper-triangular, like $\tilde{u}^1 = u^1 + \varphi_2 u^2 +  \cdots + \varphi_\ell u^\ell$ as in the proof of
Lemma~\ref{lem:endolemma}:
\begin{gather*}
\vec{B} \mapsto \begin{bmatrix}
1 & \varphi_2  & \cdots & \varphi_\ell \\
  & 1          &        &              \\
  &            & \ddots &              \\
  &            &        & 1
\end{bmatrix}
\vec{B}
\begin{bmatrix}
1 &-\varphi_2  & \cdots &-\varphi_\ell &     \\
  & 1          &        &              &  \\
  &            & \ddots &              &  \\
  &            &        & 1            &  \\
  &            &        &              &I_{n-\ell}
\end{bmatrix}.
\label{bigBtr}
\end{gather*}
The goal of this article is to understand involutivity in terms of~\eqref{bigB}.  Lemma~\ref{lem:invend} and Theorem~\ref{thm:gnf+}
accomplish this by providing a complete
version of Guillemin normal form of the type suggested but not completed in
Chapter IV, \S~5 of~\cite{BCGGG}.

\section{Involutivity}

\begin{Lemma}[linear involutivity criteria]\label{lem:invend}
Let $A$ denote a tableau given in a generic basis by
symbol relations~\eqref{eq:symrels}, as in Fig.~{\rm \ref{figtab}}.
If~$A$ is involutive, then~$A$ admits a basis of~$W$ in which it is
endovolutive, as in Fig.~{\rm \ref{figB}}.
\end{Lemma}

\begin{Theorem}[quadratic involutivity criteria]\label{thm:gnf+}
Let $A$ denote an endovolutive tableau given in a generic basis by
symbol relations~\eqref{eq:symrels}, as in Fig.~{\rm \ref{figtab}}.
The tableau~$A$ is involutive if and only if
for all~$b$, all $\lambda < i < j$ and $\lambda \leq \mu < j$, and all $a> s_i$, we have
\[\big(B^\lambda_iB^\mu_j - B^\lambda_jB^\mu_i\big)^a_b=0.\]
In particular, this implies $B(u^\lambda)(v)$ is an endomorphism of $\Wu^-_\lambda$ such
that
for all $v, \tilde{v} \in Y$,
\[
\big[B\big(u^\lambda\big)(v), B\big(u^\lambda\big)(\tilde{v})\big]=0.
\]
\end{Theorem}

A proof of Lemma~\ref{lem:invend} appears on pages 145--147 (pages~126--127 in
the online version) of~\cite{BCGGG} and it is implicit in Section~1.2 of~\cite{Yang1987}.
The proof of Theorem~\ref{thm:gnf+} is a lengthy inductive argument, each step
of which is modeled on the approach of~\cite{Yang1987}.

\begin{proof}
Choose generic bases for $W$ and $V^*$.  As a subspace, $A \subset W\otimes V^*$
is def\/ined by a~minimal set of equations
\begin{gather*}
\pi^a_i = \sum_{b\leq s_\lambda} B^{a,\lambda}_{i,b} \pi^b_\lambda,\qquad
\forall \, a > s_i,
\tag{\ref{eq:symrels} \emph{bis}}
\end{gather*}
where $B^{a,\lambda}_{i,b}=0$ unless $\lambda \leq i$ and $b \leq s_\lambda$.
Moreover, we assume that the basis of $W$ is endovolutive, so that
$B^{a,\lambda}_{i,b} =0$ if $a > s_\lambda$.

Let $\{ z^a_i,\; a \leq s_i\}$ be a basis for the abstract vector space $A^*$.
Def\/ine a monomorphism $A \to W \otimes V^*$ by
\begin{gather}
\begin{split}
\pi^a_1 &=  B^{a,1}_{1,b} z^b_1,\\
\pi^a_2 &=  B^{a,1}_{2,b} z^b_1 +
            B^{a,2}_{2,b} z^b_2,\\
\pi^a_3 &=  B^{a,1}_{3,b} z^b_1 +
            B^{a,2}_{3,b} z^b_2 +
            B^{a,3}_{3,b} z^b_3, \\
&\vdots \phantom{\sum_{b \leq s_1} B^{a,1}_{n,b} z^b_1 +
            B^{a,2}_{n,b} z^b_2 +
            \cdots + }\ddots\\
 \pi^a_n &= B^{a,1}_{n,b} z^b_1 +
            B^{a,2}_{n,b} z^b_2 +
            \cdots \quad \cdots +
            B^{a,n}_{n,b} z^b_n.
\end{split}
\label{eq:ztopi}
\end{gather}
We set $B^{a,\lambda}_{\lambda,b} = \delta^a_b$ for $a \leq s_\lambda$,
so that~\eqref{eq:ztopi} satisf\/ies~\eqref{eq:symrels}.

The prolongation $A^{(1)} \subset A \otimes V^*$ is given by coef\/f\/icients
$\{Z^a_{i,j},\; a \leq s_i\}$ with the ``contact'' system taking the form
$z^a_i=Z^a_{i,j}u^j$.  By Cartan's test, the condition of involutivity means that exactly
$s_1+2s_2+\cdots+ns_n$ of these coef\/f\/icients are independent functions on
$A^{(1)}$.

This proof is based on Section~1.1 of~\cite{Yang1987}, where it is shown that a
tableau is involutive precisely when, in generic bases, the usual 2-form condition
\begin{gather*}
0 \equiv \pi^a_1 \wedge u^1 + \pi^a_2 \wedge u^2 +
\cdots + \pi^a_n\wedge u^n
\label{2form}
\end{gather*}
is equivalent to the sequence of conditions
\begin{gather}
\begin{split}
0 &\equiv \pi^a_1 \wedge u^1 \ \mod  u^2, \ldots, u^n, \\
0 &\equiv \pi^a_1 \wedge u^1 + \pi^a_2 \wedge u^2 \ \mod  u^3, \ldots, u^n, \\
  &\phantom{=}\vdots \phantom{\pi^a_1 \wedge u^1 + \pi^a_2 \wedge u^2 + \pi^a_3 \wedge u^3}\ddots \\
0 &\equiv \pi^a_1 \wedge u^1 + \pi^a_2 \wedge u^2 + \cdots + \pi^a_k
\wedge u^k  \mod  u^{k+1},\ \ldots, u^n, \\
  &\phantom{=}\vdots \phantom{\pi^a_1 \wedge u^1 + \pi^a_2 \wedge u^2 + \cdots + \pi^a_i \wedge u^i,}\ddots
      \\
0 &\equiv  \pi^a_1 \wedge u^1 + \pi^a_2 \wedge u^2 +
 \cdots\qquad\cdots + \pi^a_n\wedge u^n.
\end{split}
\label{2formseq}
\end{gather}
The following argument shows that this sequence of conditions forces $Z^a_{i,j}$ with $a
\leq s_i$ and $j \leq i$ to be a complete set of independent generators of
$A^{(1)}$, providing it with a dimension of $s_1 + 2s_2 + 3s_3 + \cdots + n
s_n$.   To emphasize these terms in the following computations, we
\underline{underline} them.

The proof proceeds by induction, sequentially verifying that each term of
each row of \eqref{2formseq} yields a condition of the desired form.  One may
interpret this as induction over~$n$ of Cartan's test for all tableaux
of size $r \times n$.

Row~1 of \eqref{2formseq} prolongs to
\begin{gather*}
0 = \pi^a_1 \wedge u^1
  = B^{a,1}_{1,b} z^b_1 \wedge u^1
  = B^{a,1}_{1,b} Z^b_{1,i} u^i \wedge u^1
  \equiv
     B^{a,1}_{1,b} \Z{1}{1}{b} u^1 \wedge u^1 \  \mod u^2, \ldots, u^n,
\end{gather*}
which is trivial.  \emph{Every tableau with $n=1$ is
involutive.}   The $\Z{1}{1}{b}$ terms account for $s_1$ generators
of~$A^{(1)}$.

Row 2 of \eqref{2formseq} prolongs to
\begin{gather*}
0 = \pi^a_1 \wedge u^1 + \pi^a_2 \wedge u^2
  =
    \big( B^{a,1}_{1,b} z^b_1 \big)\wedge u^1
   +\big( B^{a,1}_{2,b} z^b_1
         + B^{a,2}_{2,b} z^b_2 \big)\wedge u^2\\
\hphantom{0}{} =
    \big( B^{a,1}_{1,b} Z^b_{1,i} \big)u^i\wedge u^1
   +\big( B^{a,1}_{2,b} Z^b_{1,i}
         + B^{a,2}_{2,b} Z^b_{2,i} \big)u^i\wedge u^2\\
\hphantom{0}{}
\equiv
    \big( B^{a,1}_{2,b} \Z{1}{1}{b}
         + B^{a,2}_{2,b} \Z{2}{1}{b}
         - B^{a,1}_{1,b} \Z{1}{2}{b}
         \big)u^1\wedge u^2 \ \mod u^3, \ldots, u^n.
\end{gather*}
Endovolutivity implies this is trivial when projected to $\Wu^+_1$ (meaning
``for $a > s_1$''). But, when considering the projection to $\Wu^-_1$, we see
the condition
\begin{gather*}
Z^a_{1,2} =  B^{a,1}_{2,b}\Z{1}{1}{b} +
B^{a,2}_{2,b}\Z{2}{1}{b},\qquad \forall\, a \leq s_1.
\end{gather*}
The $\Z{2}{1}{b}$ and $\Z{2}{2}{b}$ terms account for
$2s_2$ new generators of $A^{(1)}$.  So far, there is no quadratic condition on
$B^\lambda_i$; therefore \emph{all endovolutive tableaux with $n=2$ are
involutive.}

\textbf{Notation!}  It is clear we must confront a proliferation of indices
$a,b,c,\ldots$ covering $W$.  Henceforth, we suppress these indices and work
directly on $W$-valued objects.  Instead of saying ``$\forall\,
a\leq s_i$'' and ``$\forall\, a > s_i$,'' we say ``on $\Wu^-_i$'' and ``on
$\Wu^+_i$,'' respectively.  Note that this always refers to
projection on the \emph{range} of the expression, not a restriction of its
domain; by our def\/inition of $B^\lambda_i$, we may assume the domain is always~$W$.

Row 3 of \eqref{2formseq} prolongs to
\begin{gather*}
0 =
\pi_1 \wedge u^1 + \pi_2 \wedge u^2 + \pi_3 \wedge u^3\\
\hphantom{0}{} =
    \big( B^{1}_{1} z_1 \big)\wedge u^1
   +\big( B^{1}_{2} z_1
         + B^{2}_{2} z_2 \big)\wedge u^2
   +\big( B^{1}_{3} z_1
         + B^{2}_{3} z_2
         + B^{3}_{3} z_3 \big)\wedge u^3\\
\hphantom{0}{}
=
    \big( B^{1}_{1} Z_{1,i} \big)u^i \wedge u^1
   +\big( B^{1}_{2} Z_{1,i}
         + B^{2}_{2} Z_{2,i} \big)u^i \wedge u^2
   +\big( B^{1}_{3} Z_{1,i}
         + B^{2}_{3} Z_{2,i}
         + B^{3}_{3} Z_{3,i}\big)u^i \wedge u^3\\
\hphantom{0}{}
\equiv
\big(
     B^{1}_{2} \Z{1}{1}{}
   + B^{2}_{2} \Z{2}{1}{}
   - B^{1}_{1} \Z{1}{2}{}\big)u^1\wedge u^2
+ \big(
     B^{1}_{3} \Z{1}{1}{}
   + B^{2}_{3} \Z{2}{1}{}
   + B^{3}_{3} \Z{3}{1}{}
   - B^{1}_{1} \Z{1}{3}{}
   \big)u^1 \wedge u^3\\
\hphantom{0=}{}
+\big(
     B^{2}_{3} \Z{2}{2}{}
   + B^{1}_{3} \Z{1}{2}{}
   + B^{3}_{3} \Z{3}{2}{}
   - B^{1}_{2} \Z{1}{3}{}
   - B^{2}_{2} \Z{2}{3}{}  \big)u^2\wedge u^3  \
    \mod u^4,\ldots,u^n.
\end{gather*}
For this to vanish each component $u^i \wedge u^j$ with $i < j$ must vanish separately.
The $u^1\wedge u^2$ term repeats conditions already seen in row~2, namely
\begin{gather*}
\Z{1}{2}{}  =  B^{i}_{2}\underline{Z_{i,1}} \qquad \text{on $\Wu^-_1$}.
\end{gather*}
The $u^1 \wedge u^3$ term is similar,
\begin{gather*}
\Z{1}{3}{}  =  B^{i}_{3}\underline{Z_{i,1}} \qquad \text{on $\Wu^-_1$}.
\end{gather*}
The $u^2\wedge u^3$ term is more interesting, because it requires expansion
using the previous relations:
\begin{gather}
\Z{2}{3}{}
=
 B^{1}_{3}\Z{1}{2}{} + B^{2}_{3}\Z{2}{2}{} +
 B^{3}_{3}\Z{3}{2}{} - B^{1}_{2}\Z{1}{3}{}\nonumber\\
\hphantom{\Z{2}{3}{}}{}
=
 B^{1}_{3}\big( B^{1}_{2}\Z{1}{1}{} +
 B^{2}_{2}\Z{2}{1}{}\big) + B^{2}_{3}\Z{2}{2}{} +
 B^{3}_{3}\Z{3}{2}{}
 - B^{1}_{2}\big( B^{1}_{3}\Z{1}{1}{}
 + B^{2}_{3}\Z{2}{1}{} + B^{3}_{3}\Z{3}{1}{}\big)\nonumber\\
\hphantom{\Z{2}{3}{}}{}
=
 \big(B^{1}_{3}B^{1}_{2} - B^{1}_{2}B^{1}_{3}\big)\Z{1}{1}{}
+\big(B^{1}_{3}B^{2}_{2} - B^{1}_{2}B^{2}_{3}\big)\Z{2}{1}{}\nonumber\\
\hphantom{\Z{2}{3}{}=}{}
+\big(B^{1}_{3}B^{3}_{2} - B^{1}_{2}B^{3}_{3}\big)\Z{3}{1}{}
+B^{2}_{3}\Z{2}{2}{} + B^{3}_{3}\Z{3}{2}{} \qquad \text{on
$\Wu^-_1$.}
\label{eq:z23}
\end{gather}
On $\Wu^-_2$ equation~\eqref{eq:z23} merely shows how
$\Z{2}{3}{}$ depends on the previous coordinates on $A^{(1)}$.  Note that the
$\Z{3}{1}{}$, $\Z{3}{2}{}$ and $\Z{3}{3}{}$ terms
contribute another $3s_3$ generators of $A^{(1)}$.  Cartan's test fails if and
only if other relations appear among the generators $\underline{Z_{i,j}}$, $i
\geq j$.
However, on $\Wu^+_2$, many of the terms vanish by def\/inition, and the rest
impose a new \emph{quadratic} condition:
\begin{gather}
0
= \big(B^{1}_{3}B^{1}_{2} - B^{1}_{2}B^{1}_{3}\big)\Z{1}{1}{}
+\big(B^{1}_{3}B^{2}_{2} - B^{1}_{2}B^{2}_{3}\big)\Z{2}{1}{}
+\big(B^{1}_{3}B^{3}_{2} - B^{1}_{2}B^{3}_{3}\big)\Z{3}{1}{} \qquad
\text{on $\Wu^+_2$.}
\label{eq:z23q}
\end{gather}
Therefore \emph{an endovolutive tableau with $n=3$ is involutive if and only if
each term of \eqref{eq:z23q} holds on $\Wu^+_2$.}

It is useful to see another case, where things become more interesting.

Row 4 of \eqref{2formseq} prolongs to
\begin{gather*}
0 =
\pi_1 \wedge u^1 + \pi_2 \wedge u^2 + \pi_3 \wedge u^3 + \pi_4\wedge u^4\\
\hphantom{0}{} =
    \big( B^{1}_{1} z_1 \big)\wedge u^1
   +\big( B^{1}_{2} z_1
         + B^{2}_{2} z_2 \big)\wedge u^2
   +\big( B^{1}_{3} z_1
         + B^{2}_{3} z_2
         + B^{3}_{3} z_3 \big)\wedge u^3 \\
\hphantom{0=}{}
   +\big( B^{1}_{4} z_1
         + B^{2}_{4} z_2
         + B^{3}_{4} z_3
         + B^{4}_{4} z_4 \big)\wedge u^4\\
\hphantom{0}{} =
    \big( B^{1}_{1} Z_{1,i} \big)u^i \wedge u^1
   +\big( B^{1}_{2} Z_{1,i}
         + B^{2}_{2} Z_{2,i} \big)u^i \wedge u^2
   +\big( B^{1}_{3} Z_{1,i}
         + B^{2}_{3} Z_{2,i}
         + B^{3}_{3} Z_{3,i} \big)u^i \wedge u^3\\
\hphantom{0=}{}
   +\big( B^{1}_{4} Z_{1,i}
         + B^{2}_{4} Z_{2,i}
         + B^{3}_{4} Z_{3,i}
         + B^{4}_{4} Z_{4,i} \big)u^i \wedge u^4.
\end{gather*}
After expanding these terms, modulo $u^5, \ldots, u^n$, the several conditions
are found.  From the $u^1 \wedge u^4$ term:
\begin{gather*}
\Z{1}{4}{} =
 B^{1}_{4}\Z{1}{1}{} + B^{2}_{4}\Z{2}{1}{} +
 B^{3}_{4}\Z{3}{1}{} + B^{4}_{4}\Z{4}{1}{}.
\end{gather*}
This term imposes no quadratic conditions.

From the $u^2 \wedge u^4$ term:
\begin{gather}
\Z{2}{4}{}
=
 B^{1}_{4}\Z{1}{2}{} + B^{2}_{4}\Z{2}{2}{} +
 B^{3}_{4}\Z{3}{2}{} + B^{4}_{4}\Z{4}{2}{} -
 B^{1}_{2}\Z{1}{4}{}\nonumber\\
\hphantom{\Z{2}{4}{}}{}
=
B^{1}_{4}\big( B^{1}_{2}\Z{1}{1}{} +
B^{2}_{2}\Z{2}{1}{}\big) + B^{2}_{4}\Z{2}{2}{} +
B^{3}_{4}\Z{3}{2}{} + B^{4}_{4}\Z{4}{2}{} \nonumber\\
\hphantom{\Z{2}{4}{}=}{}
- B^{1}_{2}\big(
B^{1}_{4}\Z{1}{1}{} + B^{2}_{4}\Z{2}{1}{} +
B^{3}_{4}\Z{3}{1}{} + B^{4}_{4}\Z{4}{1}{}\big).
\label{eq:z24}
\end{gather}
Equation \eqref{eq:z24} becomes a new quadratic condition when projected on
$\Wu^+_2$:
\begin{gather}
0
=  B^{1}_{4}\big( B^{1}_{2}\Z{1}{1}{} +
B^{2}_{2}\Z{2}{1}{}\big)
 - B^{1}_{2}\big( B^{1}_{4}\Z{1}{1}{} + B^{2}_{4}\Z{2}{1}{} +
B^{3}_{4}\Z{3}{1}{} + B^{4}_{4}\Z{4}{1}{}\big)\nonumber\\
\hphantom{0}{}
=
\big( B^{1}_{4}B^{1}_{2}-B^{1}_{2}B^{1}_{4}\big)\Z{1}{1}{}
 +
\big( B^{1}_{4}B^{2}_{2}-B^{1}_{2}B^{2}_{4}\big)\Z{2}{1}{}
  +
\big( B^{1}_{4}B^{3}_{2}-B^{1}_{2}B^{3}_{4}\big)\Z{3}{1}{}\nonumber\\
\hphantom{0=}{}
 +
\big( B^{1}_{4}B^{4}_{2}-B^{1}_{2}B^{4}_{4}\big)\Z{4}{1}{} \qquad \text{on
$\Wu^+_2$.}
\label{eq:z24q}
\end{gather}
When reading \eqref{eq:z24q}, recall that $B^\lambda_i =0$ if $i < \lambda$.

The $u^3 \wedge u^4$ term becomes
\begin{gather}
\Z{3}{4}{}
=
 B^{1}_{4}\Z{1}{3}{} + B^{2}_{4}\Z{2}{3}{} +
 B^{3}_{4}\Z{3}{3}{} + B^{4}_{4}\Z{4}{3}{} -
 B^{1}_{3}\Z{1}{4}{} - B^{2}_{3}\Z{2}{4}{}\nonumber\\
\hphantom{\Z{3}{4}{}}{}
=
 B^{1}_{4}\big( B^{1}_{3}\Z{1}{1}{} + B^{2}_{3}\Z{2}{1}{}
 + B^{3}_{3}\Z{3}{1}{}\big)
+ B^{2}_{4}\big( B^{1}_{3}\big( B^{1}_{2}\Z{1}{1}{} +
B^{2}_{2}\Z{2}{1}{}\big) + B^{2}_{3}\Z{2}{2}{} +
B^{3}_{3}\Z{3}{2}{}\nonumber\\
\hphantom{\Z{3}{4}{}=}{}
- B^{1}_{2}\big( B^{1}_{3}\Z{1}{1}{} +
B^{2}_{3}\Z{2}{1}{} + B^{3}_{3}\Z{3}{1}{}\big)\big)
 + B^{3}_{4}\Z{3}{3}{} + B^{4}_{4}\Z{4}{3}{} \nonumber\\
\hphantom{\Z{3}{4}{}=}{}
- B^{1}_{3}\big( B^{1}_{4}\Z{1}{1}{} +
B^{2}_{4}\Z{2}{1}{} + B^{3}_{4}\Z{3}{1}{} +
B^{4}_{4}\Z{4}{1}{}\big)- B^{2}_{3}\big( B^{1}_{4}\big( B^{1}_{2}\Z{1}{1}{} +
B^{2}_{2}\Z{2}{1}{}\big) + B^{2}_{4}\Z{2}{2}{}\nonumber\\
\hphantom{\Z{3}{4}{}=}{}
+
B^{3}_{4}\Z{3}{2}{} + B^{4}_{4}\Z{4}{2}{}
 -  B^{1}_{2}\big( B^{1}_{4}\Z{1}{1}{} +
B^{2}_{4}\Z{2}{1}{} + B^{3}_{4}\Z{3}{1}{} + B^{4}_{4}\Z{4}{1}{}\big)\big).
\label{eq:z34}
\end{gather}
Equation \eqref{eq:z34} becomes a new quadratic condition when projected on
$\Wu^+_3$:
\begin{gather*}
0 =
B^{1}_{4}\big( B^{1}_{3}\Z{1}{1}{} + B^{2}_{3}\Z{2}{1}{}
 + B^{3}_{3}\Z{3}{1}{}\big)
+ B^{2}_{4}\big( B^{1}_{3}\big( B^{1}_{2}\Z{1}{1}{} +
B^{2}_{2}\Z{2}{1}{}\big) + B^{2}_{3}\Z{2}{2}{} +
B^{3}_{3}\Z{3}{2}{} \\
\hphantom{0=}{}
  - B^{1}_{2}\big( B^{1}_{3}\Z{1}{1}{} +
B^{2}_{3}\Z{2}{1}{} + B^{3}_{3}\Z{3}{1}{}\big)\big)
 - B^{1}_{3}\big( B^{1}_{4}\Z{1}{1}{} +
B^{2}_{4}\Z{2}{1}{} + B^{3}_{4}\Z{3}{1}{} +
B^{4}_{4}\Z{4}{1}{}\big)\\
\hphantom{0=}{}
- B^{2}_{3}\big( B^{1}_{4}\big( B^{1}_{2}\Z{1}{1}{} +
B^{2}_{2}\Z{2}{1}{}\big) + B^{2}_{4}\Z{2}{2}{} +
B^{3}_{4}\Z{3}{2}{} + B^{4}_{4}\Z{4}{2}{}  \\
\hphantom{0=}{}
 -  B^{1}_{2}\big( B^{1}_{4}\Z{1}{1}{} +
B^{2}_{4}\Z{2}{1}{} + B^{3}_{4}\Z{3}{1}{} + B^{4}_{4}\Z{4}{1}{}\big)\big) \qquad \text{on
$\Wu^+_3$}.
\end{gather*}
This looks like a mess, but collecting terms reveals a pattern:
\begin{gather*}
0=
\big( B^{1}_{4}B^{1}_{3} - B^{1}_{3}B^{1}_{4}
+ B^{2}_{4}\big( B^{1}_{3}B^{1}_{2} - B^{1}_{2}B^{1}_{3}\big)
- B^{2}_{3}\big( B^{1}_{4}B^{1}_{2} - B^{1}_{2}B^{1}_{4}\big) \big)
\Z{1}{1}{}\\
 \phantom{0=}{}
+ \big( B^{1}_{4}B^{2}_{3} - B^{1}_{3}B^{2}_{4}
+ B^{2}_{4}\big(B^{1}_{3}B^{2}_{2} - B^{1}_{2}B^{2}_{3}\big)
- B^{2}_{3}\big(B^{1}_{4}B^{2}_{2} - B^{1}_{2}B^{2}_{4}\big)\big)
\Z{2}{1}{}\\
 \phantom{0=}{}
+ \big( B^{2}_{4}B^{2}_{3} - B^{2}_{3}B^{2}_{4}\big)\Z{2}{2}{}
+
\big( B^{1}_{4}B^{3}_{3} - B^{1}_{3}B^{3}_{4}
+ B^{2}_{4}\big( B^{1}_{3}B^{3}_{2} - B^{1}_{2}B^{3}_{3}\big)\\
 \phantom{0=}{}
+ B^{2}_{3}\big( B^{1}_{4}B^{3}_{2} - B^{1}_{2}B^{3}_{4}\big)
\big)\Z{3}{1}{}
+ \big(B^{2}_{4}B^{3}_{3} - B^{2}_{3}B^{3}_{4}\big)\Z{3}{2}{}  \\
 \phantom{0=}{}
+
\big(B^{1}_{4}B^{4}_{3} - B^{1}_{3}B^{4}_{4} + B^{2}_{3}
\big( B^{1}_{4}B^{4}_{2} - B^{1}_{2}B^{4}_{4}\big) \big)\Z{4}{1}{}
+
\big(B^{2}_{4}B^{4}_{3} - B^{2}_{3}B^{4}_{4}\big)\Z{4}{2}{} \qquad \text{on
$\Wu^+_3$}.
\end{gather*}
Notice that, using the quadratic relations already discovered for $n=2$ and
$n=3$, we obtain
\begin{gather}
0=
  \big(B^{1}_{4}B^{1}_{3} - B^{1}_{3}B^{1}_{4}\big)\Z{1}{1}{}
+ \big(B^{1}_{4}B^{2}_{3} - B^{1}_{3}B^{2}_{4}\big)\Z{2}{1}{}
+ \big(B^{2}_{4}B^{2}_{3} - B^{2}_{3}B^{2}_{4}\big)\Z{2}{2}{}
\nonumber\\
\phantom{0=}{}
+ \big(B^{1}_{4}B^{3}_{3} - B^{1}_{3}B^{3}_{4}\big)\Z{3}{1}{}+ \big(B^{2}_{4}B^{3}_{3} - B^{2}_{3}B^{3}_{4}\big)\Z{3}{2}{}
+ \big(B^{1}_{4}B^{4}_{3} - B^{1}_{3}B^{4}_{4}\big)\Z{4}{1}{} \nonumber\\
\phantom{0=}{}
+ \big(B^{2}_{4}B^{4}_{3} - B^{2}_{3}B^{4}_{4}\big)\Z{4}{2}{} \qquad \text{on
$\Wu^+_3$}.
\label{eq:z34q}
\end{gather}
The generators $\Z{4}{1}{}$, $\Z{4}{2}{}$, $\Z{4}{3}{}$, and $\Z{4}{4}{}$ provide
$A^{(1)}$ with another~$4s_4$ dimensions.
Cartan's test fails if and
only if other relations appear among the generators~$\underline{Z_{i,j}}$, $i
\geq j$.
 \emph{An endovolutive tableau with $n=4$ is involutive if and only if
each term of \eqref{eq:z23q} and \eqref{eq:z24q} holds on $\Wu^+_2$ and each term of~\eqref{eq:z34q} holds on~$\Wu^+_3$.}

\textbf{Inductive hypothesis.}
Fix $k$ and $l<k$.
Assume for induction that the following are equivalent
\begin{enumerate}\itemsep=0pt
\item The f\/irst $k-1$ rows of the 2-form condition \eqref{2formseq} are
satisf\/ied.
\item The $s_1 + 2s_2 + \cdots + (k-1)s_{k-1}$ elements $\underline{Z^a_{j,i}}$
of $A^{(1)}$ with $a \leq s_j$ and $i \leq j < k$ are independent.
\end{enumerate}

More precisely, examining each row of \eqref{2formseq} in detail, assume for induction that the following are equivalent
\begin{enumerate}\itemsep=0pt
\item The f\/irst $k-1$ rows of the 2-form condition
\eqref{2formseq} are satisf\/ied, and the $u^i \wedge u^k$ terms vanish in the
$k$th row for all $i < l$.
\item For all $(j,i)  < (k,l)$ in lexicographic ordering on
pairs $\{ (j,i),\; i \leq j\}$, we have
\begin{gather}
   Z_{i,j} - \sum_{\mu=i}^jB^{\mu}_{j} \underline{Z_{\mu,i}} = \sum_{\lambda=1}^{i-1} \big( B^{\lambda}_{j}
   Z_{\lambda,i} - B^{\lambda}_{i} Z_{\lambda,j}\big) \qquad \text{on $\Wu^-_i$},
   \label{eq:ind-}
\end{gather}
and for all $\lambda <i < j$ and $\lambda \leq \mu \leq j$, we have
\begin{gather}
B^\lambda_i B^\mu_j - B^\lambda_j B^\mu_i=0 \qquad \text{on $\Wu^+_i$}.
\label{eq:ind+}
\end{gather}
\end{enumerate}

To perform the inductive step, we compute the $u^l \wedge u^k$ terms in the $k$th row of~\eqref{2formseq}:
\begin{gather}
0  =  \pi_1 \wedge u^1 + \cdots + \pi_l \wedge u^l +\cdots+ \pi_k \wedge u^k \nonumber\\
 \hphantom{0}{}  = \cdots + \left( \sum_{\lambda \leq l} B^{\lambda}_l z_{\lambda}
  \right) \wedge u^l + \cdots +  \left( \sum_{\mu\leq k} B^\mu_k z_\mu\right) \wedge u^k \nonumber\\
 \hphantom{0}{}
 \equiv  \left( \sum_{\lambda \leq l} B^{\lambda}_l Z_{\lambda,k}
  \right) u^k \wedge u^l  +  \left( \sum_{\mu\leq k} B^\mu_k
  Z_{\mu,l}\right) u^l \wedge u^k \
  \mod u^{k+1}, \ldots, u^n.
\label{eq:IJform}
\end{gather}
What follows is a tedious expansion and reduction of \eqref{eq:IJform} using~\eqref{eq:ind-} and~\eqref{eq:ind+}, along
the lines of what was performed for~\eqref{eq:z23},~\eqref{eq:z34} and~\eqref{eq:z24} above.  The goal is to expand~\eqref{eq:IJform} in terms of
the elements~$\underline{Z_{j,i}}$ for~$j \geq i$ that remain independent if
and only if Cartan's test holds.

For a little bit of sanity in the expansion that follows, we break up these
sums using the index ranges $\lambda_0 = 1, \ldots, l-1$, and $\mu_0 = l,
\ldots, k$.  Moreover, for every $p\geq 0$, we have nested index ranges $\lambda_{p+1}
= 1, \ldots, \lambda_{p}-1$ and $\mu_{p+1} = \lambda_p, \ldots, k$.

The vanishing of the $u^l\wedge u^k$ term of the $k$th row of
\eqref{2formseq} is equivalent to
\begin{gather}
0 =
  B^{\lambda_0}_k Z_{\lambda_0,l} + B^{\mu_0}_k \underline{Z_{\mu_0,l}}
- B^{\lambda_0}_l Z_{\lambda_0,k} - B^l_l Z_{l,k}.
\label{eq:IJterm}
\end{gather}

Rearranging terms,
\begin{gather*}
Z_{l,k} -B_k^{\mu_0} \underline{Z_{\mu_0,l}}
=
B_k^{\lambda_0} Z_{\lambda_0,l} - B_l^{\lambda_0} Z_{\lambda_0, k}
\end{gather*}
and expanding the right-hand side by the inductive hypothesis,
\begin{gather*}
=
  B^{\lambda_0}_k\big(B_l^{\lambda_1}Z_{\lambda_1,\lambda_0}-B_{\lambda_0}^{\lambda_1}Z_{\lambda_1,l}   + B_l^{\mu_1}\underline{Z_{\mu_1,\lambda_0}}\big)
-
B^{\lambda_0}_l\big(B_k^{\lambda_1}Z_{\lambda_1,\lambda_0}-B_{\lambda_0}^{\lambda_1}Z_{\lambda_1,k}+B_k^{\mu_1}\underline{Z_{\mu_1,\lambda_0}}\big)
\end{gather*}
and expanding again by inductive hypothesis,
\begin{gather*}
=
\big(
B_{k}^{\lambda_0}B_{l}^{\mu_1} - B_{l}^{\lambda_0}B_{k}^{\mu_1}
\big)
\underline{Z_{\mu_1,\lambda_0}}
+  B^{\lambda_0}_k\big(B_l^{\lambda_1}\big(B_{\lambda_0}^{\lambda_2}Z_{\lambda_2,\lambda_1}-B_{\lambda_1}^{\lambda_2}Z_{\lambda_2,\lambda_0}
  + B_{\lambda_0}^{\mu_2}\underline{Z_{\mu_2,\lambda_1}}\big)\\
\hphantom{=}{}
                 -
B_{\lambda_0}^{\lambda_1}\big(B_l^{\lambda_2}Z_{\lambda_2,\lambda_1}-B_{\lambda_1}^{\lambda_2}Z_{\lambda_2,l}
  + B_{l}^{\mu_2}\underline{Z_{\mu_2,\lambda_1}}
                 \big)\big) \\
\hphantom{=}{}
  -B^{\lambda_0}_l\big(B_k^{\lambda_1}\big(B_{\lambda_0}^{\lambda_2}Z_{\lambda_2,\lambda_1}-B_{\lambda_1}^{\lambda_2}Z_{\lambda_2,\lambda_0}
  + B_{\lambda_0}^{\mu_2}\underline{Z_{\mu_2,\lambda_1}}\big)\\
\hphantom{=}{}
                 -
B_{\lambda_0}^{\lambda_1}\big(B_k^{\lambda_2}Z_{\lambda_2,\lambda_1}-B_{\lambda_1}^{\lambda_2}Z_{\lambda_2,k}
  + B_{k}^{\mu_2}\underline{Z_{\mu_2,\lambda_1}}
                 \big)\big)
\end{gather*}
and rearranging,
\begin{gather}
=
\big(
B_{k}^{\lambda_0}B_{l}^{\mu_1} - B_{l}^{\lambda_0}B_{k}^{\mu_1}
\big)
\underline{Z_{\mu_1,\lambda_0}}\nonumber\\
\hphantom{=}{}
+
\big[
 B^{\lambda_0}_k\big(B_l^{\lambda_1}B_{\lambda_0}^{\lambda_2}-B_{\lambda_0}^{\lambda_1}B_l^{\lambda_2}\big)
-B^{\lambda_0}_l\big( B_k^{\lambda_1}B_{\lambda_0}^{\lambda_2} - B_{\lambda_0}^{\lambda_1}B_k^{\lambda_2}\big)\big]
Z_{\lambda_2,\lambda_1}\nonumber\\
\hphantom{=}{}
+\big[
  B^{\lambda_0}_lB_k^{\lambda_1}
- B^{\lambda_0}_kB_l^{\lambda_1}\Big]B_{\lambda_1}^{\lambda_2}Z_{\lambda_2,\lambda_0}
+ B^{\lambda_0}_kB_{\lambda_0}^{\lambda_1}B_{\lambda_1}^{\lambda_2}
Z_{\lambda_2,l}
+ B^{\lambda_0}_lB_{\lambda_0}^{\lambda_1}B_{\lambda_1}^{\lambda_2}
Z_{\lambda_2,k}\nonumber\\
\phantom{=}{}
+ \big[
B^{\lambda_0}_k\big(B^{\lambda_1}_lB^{\mu_2}_{\lambda_0}-B^{\lambda_1}_{\lambda_0}B^{\mu_2}_{l}\big)
-B^{\lambda_0}_l\big(B^{\lambda_1}_kB^{\mu_2}_{\lambda_0}-B^{\lambda_1}_{\lambda_0}B^{\mu_2}_{k}\big)
\big] \underline{Z_{\mu_2,\lambda_1}}.
\label{eq:prevexpansion}
\end{gather}
and canceling the $\underline{Z_{\mu_2,\lambda_1}}$ terms and expanding the
others by the inductive hypothesis,
\begin{gather*}
=
\big(
B_{k}^{\lambda_0}B_{l}^{\mu_1} - B_{l}^{\lambda_0}B_{k}^{\mu_1}
\big)
\underline{Z_{\mu_1,\lambda_0}}
+
\big[
 B^{\lambda_0}_k\big(B_l^{\lambda_1}B_{\lambda_0}^{\lambda_2}-B_{\lambda_0}^{\lambda_1}B_l^{\lambda_2}\big)
\\
\hphantom{=}{}
-B^{\lambda_0}_l\big( B_k^{\lambda_1}B_{\lambda_0}^{\lambda_2} -
B_{\lambda_0}^{\lambda_1}B_k^{\lambda_2}\big)\big]
\big( B_{\lambda_1}^{\lambda_3}Z_{\lambda_3,\lambda_2} -
B_{\lambda_2}^{\lambda_3} Z_{\lambda_3,\lambda_1} +
B^{\mu_3}_{\lambda_1}\underline{Z_{\mu_3,\lambda_2}} \big)
\\
\hphantom{=}{}
+\big[
  B^{\lambda_0}_lB_k^{\lambda_1}
- B^{\lambda_0}_kB_l^{\lambda_1}\big]
B_{\lambda_1}^{\lambda_2}
 \big( B_{\lambda_0}^{\lambda_3}Z_{\lambda_3,\lambda_2} -
B_{\lambda_2}^{\lambda_3} Z_{\lambda_3,\lambda_0} +
B^{\mu_3}_{\lambda_0}\underline{Z_{\mu_3,\lambda_2}} \big) \\
\phantom{=}{}
+ B^{\lambda_0}_kB_{\lambda_0}^{\lambda_1}B_{\lambda_1}^{\lambda_2}
\big( B_{l}^{\lambda_3}Z_{\lambda_3,\lambda_2} -
B_{\lambda_2}^{\lambda_3} Z_{\lambda_3,l} +
B^{\mu_3}_{l}\underline{Z_{\mu_3,\lambda_2}} \big) \\
\phantom{=}{}
- B^{\lambda_0}_lB_{\lambda_0}^{\lambda_1}B_{\lambda_1}^{\lambda_2}
\big( B_{k}^{\lambda_3}Z_{\lambda_3,\lambda_2} -
B_{\lambda_2}^{\lambda_3} Z_{\lambda_3,k} +
B^{\mu_3}_{k}\underline{Z_{\mu_3,\lambda_2}} \big)
\end{gather*}
and rearranging,
\begin{gather}
=
\big(
B_{k}^{\lambda_0}B_{l}^{\mu_1} - B_{l}^{\lambda_0}B_{k}^{\mu_1}
\big)
\underline{Z_{\mu_1,\lambda_0}}\nonumber
\\
\phantom{=}{}
+
\big[ B^{\lambda_0}_k B^{\lambda_1}_l \big(
 B^{\lambda_2}_{\lambda_0}B^{\lambda_3}_{\lambda_1}
-B^{\lambda_2}_{\lambda_1}B^{\lambda_3}_{\lambda_0}\big)
    -B^{\lambda_0}_l B^{\lambda_1}_k \big(
 B^{\lambda_2}_{\lambda_0}B^{\lambda_3}_{\lambda_1}
-B^{\lambda_2}_{\lambda_1}B^{\lambda_3}_{\lambda_0}\big)\nonumber\\
\phantom{=}{}
    +B^{\lambda_0}_lB^{\lambda_1}_{\lambda_0} \big(
 B^{\lambda_2}_k B^{\lambda_3}_{\lambda_1}
-B^{\lambda_2}_{\lambda_1} B^{\lambda_3}_k\big)
    -B^{\lambda_0}_k B^{\lambda_1}_{\lambda_0} \big(
 B^{\lambda_2}_l B^{\lambda_3}_{\lambda_1}
-B^{\lambda_2}_{\lambda_1}
B^{\lambda_3}_l\big)\big]Z_{\lambda_3,\lambda_2} \nonumber\\
\phantom{=}{}
-\big[
 B^{\lambda_0}_k\big(B_l^{\lambda_1}B_{\lambda_0}^{\lambda_2}-B_{\lambda_0}^{\lambda_1}B_l^{\lambda_2}\big)
-B^{\lambda_0}_l\big( B_k^{\lambda_1}B_{\lambda_0}^{\lambda_2} -
B_{\lambda_0}^{\lambda_1}B_k^{\lambda_2}\big)\big]B_{\lambda_2}^{\lambda_3} Z_{\lambda_3,\lambda_1}\nonumber\\
\phantom{=}{}
-\big[B^{\lambda_0}_lB_k^{\lambda_1} - B^{\lambda_0}_kB_l^{\lambda_1}\big]B_{\lambda_1}^{\lambda_2}B_{\lambda_2}^{\lambda_3}Z_{\lambda_3,\lambda_0}
-B^{\lambda_0}_kB_{\lambda_0}^{\lambda_1}B_{\lambda_1}^{\lambda_2}B_{\lambda_2}^{\lambda_3}Z_{\lambda_3,l}
+B^{\lambda_0}_lB_{\lambda_0}^{\lambda_1}B_{\lambda_1}^{\lambda_2}B_{\lambda_2}^{\lambda_3}Z_{\lambda_3,k}
\nonumber\\
\phantom{=}{}
+
\big[ B^{\lambda_0}_k B^{\lambda_1}_l \big(
 B^{\lambda_2}_{\lambda_0}B^{\mu_3}_{\lambda_1}
-B^{\lambda_2}_{\lambda_1}B^{\mu_3}_{\lambda_0}\big)
    -B^{\lambda_0}_l B^{\lambda_1}_k \big(
 B^{\lambda_2}_{\lambda_0}B^{\mu_3}_{\lambda_1}
-B^{\lambda_2}_{\lambda_1}B^{\mu_3}_{\lambda_0}\big)\nonumber\\
 \phantom{=}{}
    +B^{\lambda_0}_lB^{\lambda_1}_{\lambda_0} \big(
 B^{\lambda_2}_k B^{\mu_3}_{\lambda_1}
-B^{\lambda_2}_{\lambda_1} B^{\mu_3}_k\big)
    -B^{\lambda_0}_k B^{\lambda_1}_{\lambda_0} \big(
 B^{\lambda_2}_l B^{\mu_3}_{\lambda_1}
-B^{\lambda_2}_{\lambda_1}
B^{\mu_3}_l\big)\big]\underline{Z_{\mu_3,\lambda_2}}.
\label{eq:lastexpansion}
\end{gather}
The $\underline{Z_{\mu_3,\lambda_2}}$ terms cancel by the inductive
hypothesis.

Comparing \eqref{eq:lastexpansion} to \eqref{eq:prevexpansion}, it is apparent
that this pattern continues as we expand by the inductive hypothesis;
in particular, notice that the upper indices on $Z_{\lambda_{p}, \lambda_{q}}$
or $\underline{Z_{\mu_{p},\lambda_{p-1}}}$
always appear as $\lambda_0, \lambda_1, \ldots, \lambda_p$ (or $\mu_p$), while the lower indices
vary through signed permutations of $(l,k,\lambda_0, \ldots,
\lambda_q, \ldots, \lambda_p)$ that end in  $\lambda_p,\lambda_q$.
Because these indices satisfy $1 \leq \lambda_p < \lambda_{p-1} < \dots < \lambda_0 < l$, eventually
every $Z_{\lambda_p, \lambda_q}$ term will reduce by repeated application of
\eqref{eq:ind-} to terms of the
form $\underline{Z_{\mu_p,\lambda_{p-1}}}$.  Therefore,
by pairing the lower-index permutations by transposition in the third-to-last
and fourth-to-last slots, the $\underline{Z_{\mu_{p},\lambda_{p-1}}}$ terms always appear as
\begin{gather*}
\cdots
\big(B^{\lambda_{p-1}}_iB^{\mu_{p}}_j
- B^{\lambda_{p-1}}_jB^{\mu_{p}}_i\big)\underline{Z_{\mu_{p},\lambda_{p-1}}} \qquad
\text{with $\mu_{p} \geq \lambda_{p-1},\quad (j,i) < (k,l)$},
\end{gather*} which vanishes by the
inductive hypothesis.

Therefore, equation \eqref{eq:IJterm} reduces by induction to
\begin{gather*}
Z_{l,k}  -B_k^{\mu_0} \underline{Z_{\mu_0,l}} =
\big(B_{k}^{\lambda_0}B_{l}^{\mu_1} - B_{l}^{\lambda_0}B_{k}^{\mu_1}\big)
\underline{Z_{\mu_1,\lambda_0}}.
\end{gather*}
On $\Wu^+_l$, the left-hand side vanishes, so the independence of the
$s_1 + 2s_2 + \cdots + k s_k$ elements
$\underline{Z_{\mu_1,\lambda_0}}$ required by Cartan's test is equivalent to
the condition
\begin{gather*}
\big(B_{l}^{\lambda_0}B_{k}^{\mu_1} - B_{k}^{\lambda_0}B_{l}^{\mu_1}\big)
=0
\end{gather*}
projected to $\Wu^+_l$ for all $\lambda_0 \leq \mu_1  \leq k$ and $\lambda_0 < l$.

This concludes the proof of Theorem~\ref{thm:gnf+}.
\end{proof}

Another important subspace is the part of $A$ that is rank-one in $W
\otimes U^*$; that is, consider
\begin{gather*}
\Au^1(\varphi) = A \cap \big\{ z \otimes \varphi  + J
\text{ for some $z \in W$, $J \in W \otimes Y^*$}  \big\}.
\end{gather*}
The image of $\Au^1(\varphi)$ under the projection $(W \otimes V^*)
\to W \otimes U^*$ is comprised of rank-one homomorphisms, so
the projection $\Au^1(\varphi) \to W$ is well-def\/ined, with image
\begin{gather*}
\Wu^1(\varphi) = \big\{ z \in W\colon \; z\otimes \varphi +J \in
A \text{ for some $J \in W \otimes Y^*$} \big\}.
\label{W1e}
\end{gather*}
The spaces $\Wu^-(\varphi)$ and $\Wu^1(\varphi)$ are distinct, but their relationship is clear:

\begin{Lemma}\label{lem:Wup}
Suppose that $A$ is an endovolutive tableau.
For any $\lambda$,
\[
\Wu^1(u^\lambda) = \big\{
z \in \Wu^-_\lambda\colon \; B^\lambda_\mu z = \delta^\lambda_\mu z \; \forall \,\mu
\leq \ell\big\}.
\]
More generally, for any $\varphi \in U^*$,
\[
\Wu^1(\varphi) = \left\{
z \in \Wu^-(\varphi)\colon \;
\left( \sum_\lambda \varphi_\lambda B^\lambda_\mu - \varphi_\mu I  \right) z = 0\;
 \forall \, \mu
\leq \ell
\right\}.
\]
\end{Lemma}

There is nothing to prove; Lemma~\ref{lem:Wup} merely states the condition that the $W \otimes U^*$ part of
$\pi=B(\varphi)(\cdot)z$ is rank-one.
\begin{Corollary}
Suppose that $A$ is an endovolutive tableau.
For almost all $\varphi$, $\dim \Wu^1(\varphi) = s_\ell$.
\label{cor:dimW1}
\end{Corollary}

\begin{proof}
For almost all $\varphi \in U^*$, we have $\min \{ i\colon \varphi_i \neq 0 \} = 1$
and $\max \{ i\colon  \varphi_i \neq 0 \} = \ell$.
So, $\Wu^-(\varphi) = \Wu^-_1$ has dimension $s_1$. For each $\mu=1,\ldots,\ell$
the condition $ ( \varphi_\lambda B^\lambda_\mu - \varphi_\mu
I   ) z = 0$ has no $\Wu^-_\mu$ component, but the
$\Wu^+_\mu$ component has rank up to $s_{\mu-1} - s_{\mu}$.  The rank falls if
and only if
$\varphi_\mu$ is an eigenvalue of $\sum\limits_{\lambda < \mu} \varphi_\lambda
B^\lambda_\mu$.  By f\/ixing $\varphi_1$, $\varphi_2$,
$\ldots$, $\varphi_{\mu-1}$ and varying~$\varphi_\mu$, we can see that this
condition achieves its maximum rank for a Zariski-open set of values of
$\varphi$. Then
$\dim \Wu^1(\varphi) = s_1 - (s_1-s_2) -  \cdots - (s_{\ell-1} -
s_{\ell}) = s_\ell$.
\end{proof}

Note that, unlike with $\Wu^-(\varphi)$, the def\/inition of $\Wu^1(\varphi)$
does not rely on the basis; it requires only a splitting of $0 \to U^* \to V^*
\to Y \to 0$ to decide where $J$ takes values.  For this
reason, it is the space studied in the homological references.
The space $\Wu^1(\varphi)$ is the subject of Theorem~\ref{thm:gnf1},
which is known colloquially as ``Guillemin normal form''.
It is called\footnote{In these references, the domain is restricted to $v \in Y$, but that limitation is artif\/icial.} Lemma~4.1 in~\cite{Guillemin1968} and Proposition~6.3 in Chapter~VIII of~\cite{BCGGG}.
Theorem~\ref{thm:gnf1} is crucial to the study of partial
dif\/ferential equations and exterior dif\/ferential systems because it reveals the intimate
relationship between involutivity, overdetermined Cauchy initial-value
problems, and the characteristic variety~\cite{Guillemin1970}.

\begin{Theorem}[Guillemin]\label{thm:gnf1}
Suppose that $A$ is involutive.
For every $\varphi \in U^*$ and $v \in V$, the restricted homomorphism
$B(\varphi)(v)|_{\Wu^1(\varphi)}$ is an endomorphism of $\Wu^1(\varphi)$.
Moreover, for all \mbox{$v, \tilde{v} \in V$},
\begin{gather} \big[
B(\varphi)(v),
B(\varphi)(\tilde{v})\big]\big|_{\Wu^1(\varphi)}
=0.
\label{eq:gnf1}
\end{gather}
\end{Theorem}

Guillemin's original proof of Theorem~\ref{thm:gnf1} relies on several subtle homological
results.  Each of those results can be reproven using Theorem~\ref{thm:gnf+}
and elementary linear algebra.  For example, here is the key result:

\begin{Corollary}[Quillen, Guillemin]\label{thm:guilA}
If $A$ is involutive, then $A|_U$ is involutive, and the natural map
between prolongations $A^{(1)} \to  (A|_U )^{(1)}$ is bijective.
\end{Corollary}

\begin{proof}
This is Theorem~A in \cite{Guillemin1968}, where it
is proven with a large
diagram chase using Quillen's exactness theorem from \cite{Quillen1964}.  But, using Theorem~\ref{thm:gnf+}, this
is immediate, as the quadratic condition still holds if the range of indices
$\lambda$, $\mu$, $i$, $j$ is truncated at $\ell$ (or greater).  In particular, the generators $( z^a_\mu )_{a \leq s_\mu}$ of $A$ are preserved.
As explored in the proof of Theorem~\ref{thm:gnf+}, the contact relation $z^a_\mu = Z^{a}_{\mu,i} u^i$ gives
coordinates $Z^a_{\mu,i}$ to the prolongation $A^{(1)} \subset A \otimes V^*$,
and the $s_1 + 2s_2 + \cdots + \ell s_\ell$ independent generators are
$\underline{Z^a_{\mu,\lambda}}$ with $a \leq s_\mu$ and
$\lambda  \leq \mu$.  These generators remain independent under restriction
to $U$, too.
\end{proof}

Theorems~\ref{thm:gnf+} and~\ref{thm:gnf1} appear to be very
similar in the sense that each implies linear and quadratic
conditions on the symbol coef\/f\/icients of involutive tableaux. However, Theorem~\ref{thm:gnf+} is strictly
stronger, as it provides equivalence, whereas Theorem~\ref{thm:gnf1} is a~unidirectional implication.
For example, consider this tableau with characters
$(3 , 1, 0)$:
\begin{gather*}
\pi =
\begin{bmatrix}
\pi^1_1 & \pi^1_2 & P_{1} \pi^1_1 + P_{2} \pi^2_1 + P_{3} \pi^3_1 + Q \pi^1_2 \\
\pi^2_1 & \pi^3_1 & T_{2} \pi^2_1 + T_{3} \pi^3_1 \\
\pi^3_1 & 0 & R_{3} \pi^3_1
\end{bmatrix}
\end{gather*}
with symbol maps arranged as
\begin{alignat*}{4}
&\begin{pmatrix}
1 & 0 & 0 \\
0 & 1 & 0 \\
0 & 0 & 1
\end{pmatrix},
\qquad &&
\begin{pmatrix}
0 & 0 & 0 \\
0 & 0 & 1 \\
0 & 0 & 0
\end{pmatrix},
\qquad & &
\begin{pmatrix}
P_1 & P_2 & P_3 \\
  0 & T_2 & T_3 \\
  0 &   0 & R_3
\end{pmatrix}, & \\
& &&
\begin{pmatrix}
1 & ~ & ~  \\
~ & ~ & ~  \\
~ & ~ & ~  \\
\end{pmatrix},
\qquad &&
\begin{pmatrix}
Q & ~ & ~  \\
\phantom{P_1}  & \phantom{P_2}  & \phantom{P_3} \\
~  & ~  & ~
\end{pmatrix}.&
\end{alignat*}
This tableau is written in generic bases, and it is endovolutive.
It has $\Wu^-(u^1) = \pair{ w_1, w_2, w_3}$, $\Wu^-(u^1) = \pair{ w_1, w_2}$, and
$\Wu^1(u^1) = \pair{w_1,w_2}$.  For $\varphi = \varphi_1 u^1 + \varphi_2 u^2 \in
U^*$ with $\varphi_2 \neq 0$, it has $\Wu^1(\varphi) = \pair{w_1}$.
Moreover, for all $\varphi = \varphi_1 u^1 + \varphi_2 u^2 \in U^*$ and all $v
\in V$, the restricted map $B(\varphi)(v)|_{\Wu^1(\varphi)}$ is an
endomorphism of $\Wu^1(\varphi)$, and the commutativity condition~\eqref{eq:gnf1} holds.
However, by Theorem~\ref{thm:gnf+}, this tableau is involutive if and only if $T_2=R_3$.

Theorem~\ref{thm:gnf1} sets conditions on operators acting in dimension
$s_\ell$, whereas Theorem~\ref{thm:gnf+} sets conditions on operators acting in
dimension $s_1$. These additional conditions allow Theorem~\ref{thm:gnf+} to describe the
variety of involutive tableaux, a strict subvariety of the endovolutive tableaux satisfying \eqref{eq:gnf1}.

\section{Discussion}

Theorem~\ref{thm:gnf+} is the f\/irst step to answering a very fundamental open
question, which is expressed in footnote~7 in Chapter~IV of \cite{BCGGG}:
``What is the dimension of the space of involutive tableaux with certain f\/ixed
Cartan characters?''  Lemma~\ref{lem:invend} and Theorem~\ref{thm:gnf+}
together give an algebraic ideal whose variety contains every involutive
tableaux, each expressed in a preferred basis.  It is not immediately clear how to
account for non-Borel basis changes or in what sense the endovolutive expression
of $B^\lambda_i$ is unique for a given abstract tableau, but the recursive
nature of the proof suggests that Theorem~\ref{thm:gnf+} provides
\emph{minimal} set of generators for this ideal.

Because it is easy to program into computer algebra systems\footnote{I am
writing an open-source package called \texttt{Symbol} for Sage using this
approach.  You can help at \url{https://bitbucket.org/curieux/symbol_sage}.},
Theorem~\ref{thm:gnf+} allows us to explore the moduli of involutive tableaux
and search for interesting new families of involutive partial dif\/ferential
equations with peculiar geometric properties.  More generally, I hope it
lowers the barrier for future researchers of exterior dif\/ferential systems and Lie
pseudogroups who want to understand and \emph{apply} the profound results
of \cite{Quillen1964} and \cite{Guillemin1968}.

\subsection*{Acknowledgments}
Thanks to Deane Yang for several helpful conversations.  Thanks also to the
anonymous referees, whose suggestions improved the style and focus of this
article signif\/icantly.

\pdfbookmark[1]{References}{ref}
\LastPageEnding

\end{document}